%% file: split.tex
\documentclass[12pt]{article}

\usepackage{epsf}
\usepackage{psfrag}
\usepackage{amsmath}
\usepackage{amsthm}
\usepackage{amssymb}
\usepackage{delarray}

\newcolumntype{L}{>{$}l<{$}}
\newenvironment{Cases}{\begin{array}\{{lL}.}{\end{array}}

\renewcommand{\epsfsize}[2]{0.7\textwidth}

\theoremstyle{plain}
\newtheorem{thm}{Theorem}[section]
\newtheorem{prop}[thm]{Proposition}
\newtheorem{cor}[thm]{Corollary}

\newtheorem{intlemnp}[thm]{Lemma}

\newtheorem{intthmnp}[thm]{Theorem}
\newenvironment{thmnp}{\begin{intthmnp}}{\qed \end{intthmnp}}
\newtheorem{intpropnp}[thm]{Proposition}

\theoremstyle{definition}

\newtheorem{ex}[thm]{Example}

\theoremstyle{remark}
\newtheorem{rem}[thm]{Remark}

\def\thmref#1{Theorem~\ref{#1}}

\def\propref#1{Proposition~\ref{#1}}

\def\figref#1{Figure~\ref{#1}}

\def\secref#1{Section~\ref{#1}}

\title{Splittings of mapping tori of free group automorphisms}
\date{\today}
\author{Peter Brinkmann}

\begin{document}
\renewcommand{\thefootnote}{\null}
\maketitle%
\footnote{{\em 2000 Mathematics Subject Classification.} 37E30.}
\footnote{{\em Key words and phrases.} Hyperbolic groups, Gromov boundary,
splittings over $\mathbb{Z}$.}
\footnote{Research partially supported by a Liftoff Grant from
the Clay Mathematics Institute.}
\setcounter{footnote}{0}
\renewcommand{\thefootnote}{\arabic{footnote}}

\input{abstract}

\input{intro}

\input{splitsec}

\input{examples}

\bibliographystyle{alpha}
\bibliography{$GLOBAL/all}
\par

{\noindent \sc Department of Mathematics\\
University of Illinois at Urbana-Champaign\\
273 Altgeld Hall\\
1409 W.\ Green St.\\
Urbana, IL 61801, USA\\}
\par
{\noindent \it E-mail:} brinkman@math.uiuc.edu

\end{document}

%% file: abstract.tex
\begin{abstract}
We present necessary and sufficient conditions for the existence of a
splitting over $\mathbb{Z}$ of the mapping torus
$M_\phi=F\rtimes_\phi \mathbb{Z}$ of a free group automorphism $\phi$.
\end{abstract}

%% file: intro.tex
\section*{Introduction}

Let $\phi$ be an automorphism of a finitely generated free group $F$ of
rank at least two, and
let $M_\phi=F\rtimes_\phi\mathbb{Z}$ be the mapping torus of $\phi$.
If $F$ is generated by $x_1,\ldots,x_b$, then $M_\phi$ is presented by
\[
\langle~x_1,\ldots,x_b,t~|~tx_it^{-1}=\phi(x_i)~\rangle.
\]

An automorphism $\phi$ is called {\em hyperbolic} if $M_\phi$ is
a word-hyperbolic group, and it is called {\em atoroidal} if no
positive power of $\phi$ preserves the conjugacy class of a nontrivial
element of $F$. The following result gives us a nice way
of detecting hyperbolicity.

\begin{thm}[\cite{pbhyp}]\label{thesismain}
An automorphism $\phi: F\rightarrow F$ is hyperbolic if and only
if it is atoroidal.
\end{thm}

I.\ Kapovich \cite{ilyaendo} has proved a similar result for mapping
tori of certain injective endomorphisms of free groups.

Given a word-hyperbolic mapping torus $M_\phi$, it is natural to ask
whether it has a nontrivial JSJ decomposition \cite{jsj} (see also
\cite{scswjsj}).
This question appeared as Problem~6.7 in \cite{ilyaendo}. In particular,
it is natural to ask whether $M_\phi$ splits over $\mathbb{Z}$, i.e.,
whether $M_\phi$ can be expressed as an HNN extension or a nontrivial
free product with amalgamation over $\mathbb{Z}$.

The main result of this paper determines exactly when $M_\phi$ splits
in this fashion (\thmref{splitthm}). This answers a question of Lee Mosher.
The actual result is rather technical, but its main interest lies in
the fact that splittings over $\mathbb{Z}$ can arise in various different
and unexpected ways (contrast this with mapping tori of automorphisms
of surface groups, which never split over $\mathbb{Z}$).

Note that $M_\phi$ does not split as a nontrivial free product, i.e.,
it is one-ended. Moreover, a result of M.\ Kapovich and B.\ Kleiner
\cite[Theorem 14]{kk} implies that the Gromov boundary
$\partial_\infty M_\phi$ is the Menger curve if $M_\phi$ is hyperbolic
and does not split over $\mathbb{Z}$.

There exists an algorithm for finding the JSJ decomposition of a one-ended,
torsion-free hyperbolic group \cite{Sela95}. The results of this paper are
nonconstructive, but it is conceivable that there exists a simpler (and
faster) algorithm for finding splittings of mapping tori of free group
automorphisms.  Another open problem is the extension of the results of
this paper to the case of mapping tori of injective endomorphisms of free
groups.

\secref{splitsec} contains the statement and proof of the main result.
In \secref{exsec}, we give some examples of hyperbolic automorphisms whose
mapping tori split over $\mathbb{Z}$. One interesting consequence of these
examples is that hyperbolic automorphisms may have train track representatives
\cite{hb1} with strata of polynomial growth.

I would like to thank Lee Mosher for suggesting the problem, and
John Stallings and Ilya Kapovich for helpful discussions. The referee
made several suggestions that improved the paper. I am indebted to the
University of Osnabr\"uck as well as the University
of California at Berkeley for their hospitality, and to the Clay
Mathematics Institute for financial support.

%% file: splitsec.tex
\section{Splittings of $M_\phi$ over $\mathbb{Z}$}\label{splitsec}

Unless stated otherwise, we assume throughout this section
that $M_\phi$ splits over $\mathbb{Z}$,
either as an HNN-extension or as a free product with amalgamation. Hence,
$M_\phi$ acts (on the left) on a tree $T$ with one orbit of edges and
cyclic edge stabilizers \cite{arbres} (see \cite{topics} or \cite{scottwall}
for a brief review of the facts from Bass-Serre theory that we use here).
Let $t$ denote the stable letter of $M_\phi$. For a vertex $V$ of $T$, let
$F_V=F\cap~{\text{Stab}}~V$. Given an edge $E$ of $T$, we write $\tau(E)$
for the terminal vertex of $E$ and $\iota(E)$ for the initial vertex.

Let $\Gamma=F\setminus T$ be the quotient of $T$ under the action of $F$.
$\Gamma$ expresses $F$ as the fundamental group of a finite
graph of groups with trivial edge stabilizers. In order to see this,
suppose that $\Gamma$ is not finite, i.e., $\Gamma$ expresses $F$ as
the fundamental group of an infinite graph of groups with cyclic edge
groups. As $F$ is normal in $M_\phi$, there is an induced action of the
stable letter $t$ on $\Gamma$, and the quotient of $\Gamma$ under the
action of $\langle t \rangle$ consists of one edge. This forces $\Gamma$
to be homeomorphic to the real line, all vertex groups are conjugates
of each other, and all edge groups are conjugates of each other. Hence,
if the vertex groups are strictly larger than the edge groups, then $F$ is 
infinitely generated, and if all the vertex groups are the same as the
edge groups, then $F$ is cyclic. Since $F$ is finitely generated and has
rank at least two, both cases are impossible.

A careful analysis of the action of the stable letter $t$ on $\Gamma$ will
give us a description of all possible splittings of mapping tori over
$\mathbb{Z}$.  The computational techniques are similar in all cases.
In order to avoid redundancy, we spell out the details of the
computations in the first few cases and restrict the exposition to an
outline once all the ideas have appeared.

\subsection{HNN extension}

We assume that $M_\phi$ splits as an HNN extension
over $\mathbb{Z}$, i.e., $M_\phi=G\ast_{\mathbb{Z}}$. In particular,
$t$ acts on $\Gamma$ with one orbit of vertices and one orbit of edges.
We distinguish three subcases depending on the translation length of the
action of $t$ on $T$.

\subsubsection{Elliptic $t$-action}

\begin{figure}
\renewcommand{\epsfsize}[2]{0.8\textwidth}
\psfrag{V}{$V$}
\psfrag{v}{$v$}
\psfrag{e}{$E$}
\psfrag{te}{$tE$}
\psfrag{t2e}{$t^2E$}
\psfrag{t3e}{$\cdots$}
\psfrag{t4e}{$t^kE$}
\center{\epsfbox{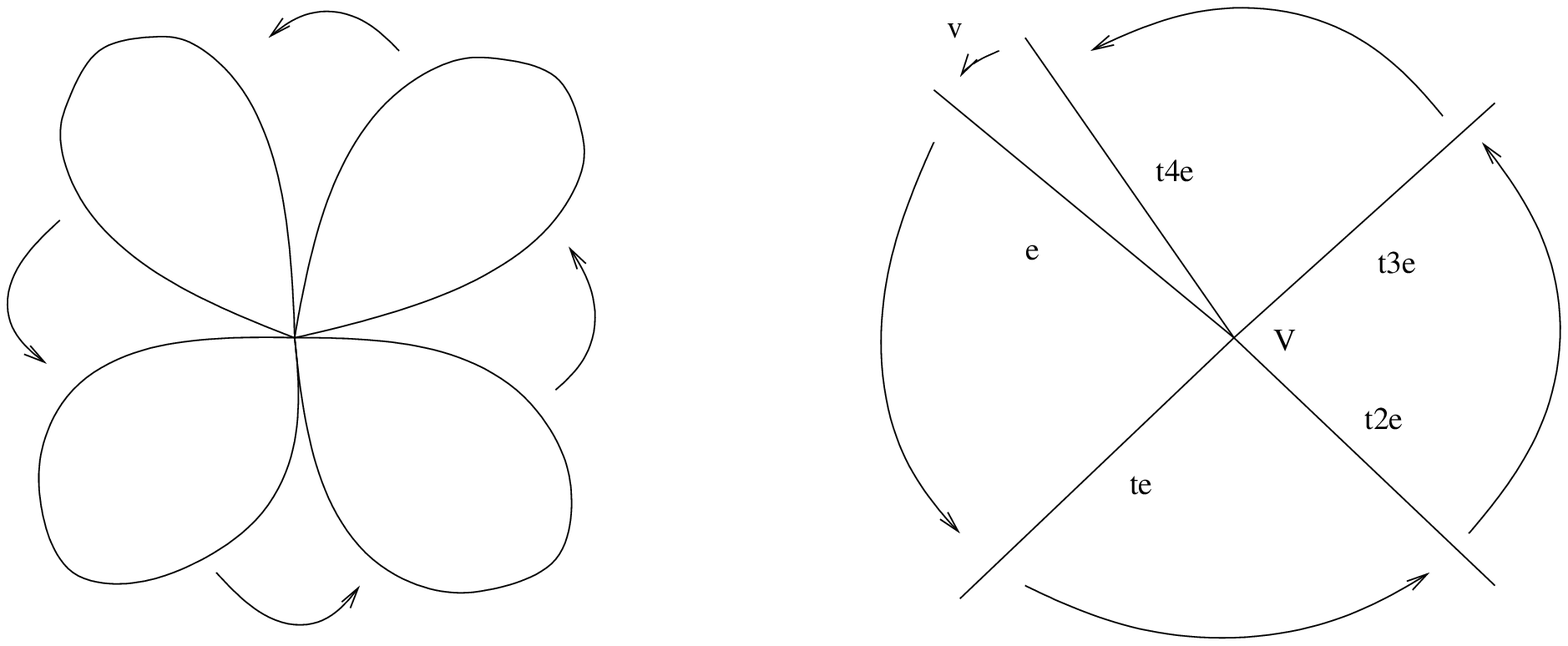}}
\caption{HNN case, elliptic $t$-action.}\label{case1pic}
\end{figure}

We consider quotients $\Gamma$ with one vertex and one orbit of $k$ edges
(\figref{case1pic}). In this case, $t$ acts on $T$ as an elliptic isometry,
and it fixes exactly one vertex $V$ of $T$. We have
\[
F=F_V\ast \langle a_0,\ldots,a_{k-1} \rangle
\]
such that
\begin{align*}
\phi(F_V)&=F_V\\
\phi(a_i)&=\begin{Cases}
a_{i+1} & if $0\leq i< k-1$;\\
wa_0v   & if $i=k-1$, for some $v,w\in F_V$.\\
\end{Cases}
\end{align*}

\begin{proof}
Clearly, if $x\in F_V$, then $V=txV=\phi(x)tV=\phi(x)V$, hence
$\phi(x)\in F_V$.

Let $E$ be some edge edge with $\iota(E)=V$ ($V$ is the vertex fixed by
$t$). Choose $a_0,\ldots,a_{k-1}\in F$ such that $a_i\tau(t^iE)=V$, and
choose $v\in F_V$ such that $vt^kE=E$ (\figref{case1pic}).
With these definitions, we have
\[
F=F_V\ast\langle a_0,\ldots,a_{k-1} \rangle.
\]

We have $ta_i\tau(t^iE)=tV$, which implies $\phi(a_i)\tau(t^{i+1}E)=V$.
Hence, for $i<k-1$, there is no loss in assuming that $\phi(a_i)=a_{i+1}$
(otherwise, we could simply modify our choice of $a_{i+1}$).
Finally, we have $\phi(a_{k-1})\tau(t^kE)=\phi(a_{k-1})v^{-1}\tau(E)=V$,
which implies $a_0^{-1}V=v\phi(a_{k-1}^{-1})V$, which in turn implies
$\phi(a_{k-1})=wa_0v$ for some $w\in F_V$.
\end{proof}

Conversely, any automorphism $\phi$ of the form listed above will
give rise to a mapping torus $M_\phi$ that splits over $\mathbb{Z}$. To
show this, we will give a sequence of Tietze transformations that exhibits
the splitting.

\begin{equation*}
\begin{split}
M_\phi &= \langle F_V,a_0,\ldots,a_{k-1},t~\vert~tF_Vt^{-1}=F_V,
ta_0t^{-1}=a_1,\ldots,ta_{k-1}t^{-1}=wa_0v\rangle \\
&\cong \langle F_V,a_0,t~\vert~txt^{-1}=\phi(x)~\forall x\in F_V,
t^ka_0t^{-k}=wa_0v \rangle \\
&\cong \langle F_V,a_0,t~\vert~txt^{-1}=\phi(x)~\forall x\in F_V,
a_0^{-1}w^{-1}t^ka_0=vt^k
\rangle \\
&\cong \langle F_V,t~\vert~txt^{-1}=\phi(x)~\forall x\in F_V
\rangle \ast_{%
\langle w^{-1}t^k\sim vt^k \rangle}.
\end{split}
\end{equation*}

\subsubsection{HNN extension: Hyperbolic $t$-action with translation distance
$1$}\label{hnndist1}

\begin{figure}
\renewcommand{\epsfsize}[2]{0.5\textwidth}
\psfrag{v0}{$V$}
\psfrag{v1}{$tV$}
\psfrag{v2}{$t^2V$}
\psfrag{v3}{$t^3V$}
\psfrag{v4}{$t^{m-1}V$}
\center{\epsfbox{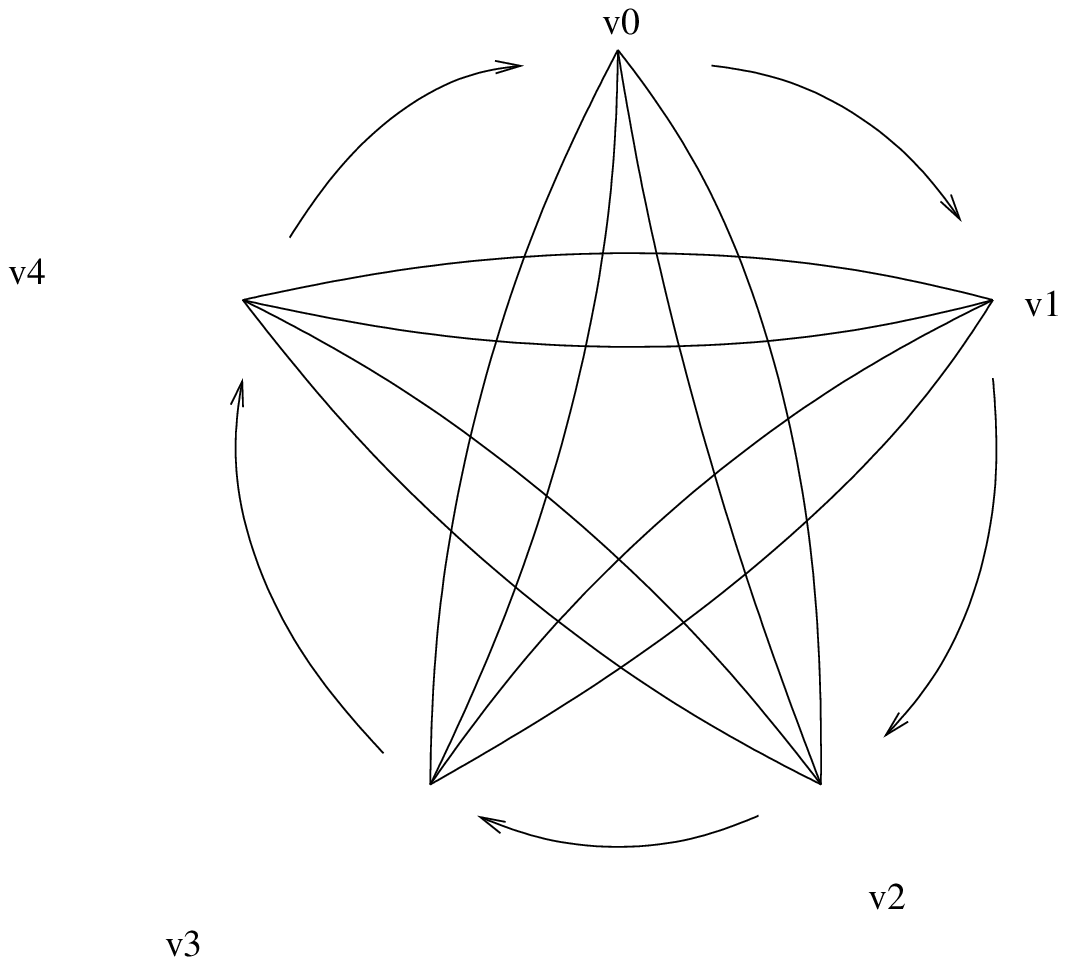}}
\caption{HNN case, hyperbolic $t$-action.}\label{case2pic}
\end{figure}

We consider the case where $\Gamma$ has one orbit of $m$ vertices
and one orbit of $km$ edges (\figref{case2pic}), and we assume
that $t$ acts on $T$ as a hyperbolic isometry with translation distance $1$.

Then there exists some vertex $V\in T$ and some edge $E\in T$ such that
$\iota(E)=V$ and $\tau(E)=tV$.  Let $v(i)$ be the remainder of $i$ under
division by $m$. We have
\[
	F=\left(\prod_{i=0}^{m-1} F_{t^iV}\right)\ast
	\langle a_{m-1},\ldots,a_{km-1} \rangle
\]
such that
\[
\phi (F_{t^iV})=\begin{Cases}
F_{t^{i+1}V} & if $0\leq i<m-1$;\\
a_{m-1}^{-1}F_{V}a_{m-1} & if $i=m-1$,\\
\end{Cases}
\]
and
\[
\phi (a_i)=\begin{Cases}
a_{i+1} & if $v(i)< m-2$;\\
a_{m-1}^{-1}a_{i+1} & if $v(i)=m-2$;\\
a_{i+1}a_{m-1}  & if $v(i)=m-1$, $i<km-1$;\\
va_{m-1} & if $i=km-1$, for some $v\in F_V$.\\
\end{Cases}
\]

\begin{figure}
\renewcommand{\epsfsize}[2]{0.9\textwidth}
\psfrag{e0}{$E_0$}
\psfrag{e1}{$E_1$}
\psfrag{e2}{$E_2$}
\psfrag{e3}{$\cdots$}
\psfrag{e4}{$E_{m-1}$}
\psfrag{tem}{$tE_{m-1}$}
\psfrag{em}{$E_m$}
\psfrag{a}{$a$}
\psfrag{V}{$V$}
\psfrag{t}{$t$}
\center{\epsfbox{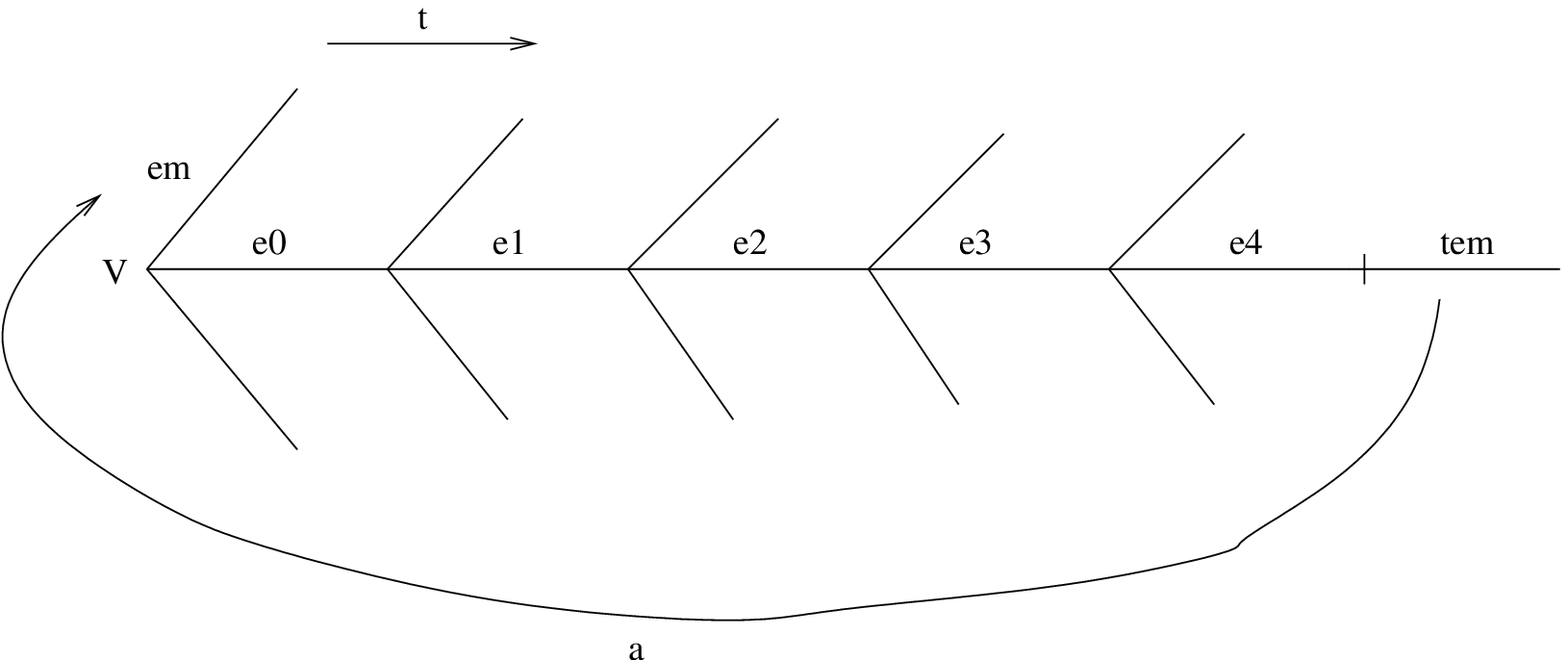}}
\caption{A lift of $\Gamma$ (case \ref{hnndist1}).}\label{case21pic}
\end{figure}

\begin{proof}

We first construct a subtree $S\subset T$ with $km$ edges that
projects to $\Gamma$. Choose $a\in F$ such that $at^mV=V$.
Let $E_0=E$, and define $E_{i+1}=tE_i$ if $v(i)\neq m-1$, and
$E_{i+1}=atE_i$ if $v(i)=m-1$. Let $S$ be the union of the
edges $E_0,\ldots,E_{km}-1$ (\figref{case21pic}, note that
the union $S'$ of $E_0,\ldots,E_{m-2}$ projects to a spanning
tree of $\Gamma$).

Now, choose $a_{m-1},\ldots,a_{km-1}\in F$ such that
$a_i\tau(E_i)=t^{v(i+1)}V$ (we choose $a_{m-1}=a$). With
these definitions, we have
\[
F=\left(\prod_{i=0}^{m-1} F_{t^iV}\right)\ast
\langle a_{m-1},\ldots,a_{km-1} \rangle,
\]
and our claim concerning the images of vertex stabilizers follows
immediately. We still need to compute the images of the generators
$a_i$.

We have $a_i\tau(E_i)=t^{v(i+1)}V$, which implies
$\phi(a_i)\tau(tE_i)=t^{v(i+1)+1}V$. If $v(i)<m-2$, we have $tE_i=E_{i+1}$
and $v(i+1)+1=v(i+2)$, so there is no loss in assuming that
$\phi(a_i)=a_{i+1}$.

If $v(i)=m-2$, then $t^{v(i+1)+1}V=t^mV=a^{-1}V$ and $tE_i=E_{i+1}$, which
implies $\phi(a_i)\tau(E_{i+1})=a^{-1}V$, so we can let
$\phi(a_i)=a^{-1}a_{i+1}$.

If $v(i)=m-1$ and $i\neq km-1$, then $tE_i=a^{-1}E_{i+1}$ and
$t^{v(i+1)+1}V=t^{v(i+2)}V$, so we have
$\phi(a_i)t\tau(E_i)=\phi(a_i)\tau(a^{-1}E_{i+1})=t^{v(i+2)}V$, hence
$\phi(a_i)a^{-1}=a_{i+1}$. This implies $\phi(a_i)=a_{i+1}a$.

Finally, there exists some $v\in F_V$ such that $vatE_{km-1}=E_0$, and
an argument similar to the previous one shows that $\phi(a_{km-1})=va$.

\end{proof}

Conversely, any automorphism $\phi$ of the form listed above will
give rise to a mapping torus $M_\phi$ that splits over $\mathbb{Z}$. To
show this, we will first find a new set of generators that will simplify
this verification. Let $b_i=a_{m+i-1}a_{m+i-2}\cdots a_{m-1}$ for
$i=0,\ldots,(k-1)m$. A simple induction shows that
\[
\phi(b_i)=\begin{Cases}
b_{i+1} & if $v(i)<m-1$ and $i\neq (k-1)m$;\\
b_0^{-1}b_{i+1} & if $v(i)=m-1$;\\
vb_{(k-1)m} & if $i=(k-1)m$.
\end{Cases}
\]

Moreover, we have
\[
\begin{split}
(a_{m-1}t^m)^k=(b_0t^m)^k=b_0\phi^m(b_0)\phi^{2m}(b_0)\cdots
\phi^{(k-1)m}(b_0)t^{km}
\\=b_0b_0^{-1}b_mb_m^{-1}b_{2m}b_{2m}^{-1}\cdots b_{(k-2)m}b_{(k-2)m}^{-1}
b_{(k-1)m}t^{km}=b_{(k-1)m}t^{km}.
\end{split}
\]
After these preliminary computations, we can list a sequence of Tietze
moves that exhibits the splitting over $M_\phi$ of $\mathbb{Z}$.

\begin{equation*}
\begin{split}
M_\phi &= \langle F_V,F_{tV},\ldots,F_{t^{m-1}V},
b_0,\ldots,b_{(k-1)m},t~\vert~\ldots \rangle \\
&\cong \langle F_V,b_0,s,t~\vert~s=b_0t^m, sF_Vs^{-1}=F_V,\ldots \rangle \\
&\cong \langle F_V,s,t~\vert~sF_Vs^{-1}=F_V,ts^kt^{-1}=vs^k \rangle\\
&\cong \langle F_V,s~\vert~sF_Vs^{-1}=F_V \rangle
\ast_{\langle s^k\sim vs^k\rangle}.\\
\end{split}
\end{equation*}
Note, in particular, that the last expression allows us to see the splitting
geometrically.

\subsubsection{HNN extension: Hyperbolic $t$-action with translation distance
greater than $1$}

The main idea of the previous argument was to find a lift $S\subset T$
of the graph $\Gamma$ with the property that a spanning tree $S'\subset S$
is contained in the axis of $t$. This particular lift of $\Gamma$ allowed
us to read off the exact appearance of the automorphism $\phi$. We will
follow the same approach here.

As in \ref{hnndist1}, $\Gamma$ has $m$ vertices and $km$ edges.
Now, however, given a vertex $V\in T$ and an edge $E$ with $\iota(E)=V$,
we have $\tau(E)=xt^rV$ for some $x\in F$ and $r>1$. Note that $m$ and
$r$ are coprime because $\Gamma$ is connected and has only one orbit
of edges. Let $s$ be the smallest positive number such that
$rs\equiv \pm 1\mod m$. Then the translation distance of $t$ equals
$d=\min\{s,m-s\}$. By inverting $t$ if necessary, we may assume that
$d=s$.

Label the edges of $\Gamma$ by letting $E_i=t^iE$ for $i=0,\ldots,km-1$.
Similarly, let $V_j=t^jV$ for $j=0,\ldots,m-1$. As before, let $v(i)$
be the remainder of $i$ under division by $m$. Note that
$\iota(E_j)=V_{v(rj)}$ and $\tau(E_j)=V_{v\left(r(j+1)\right)}$.
Hence, the edges $E_{r-1},E_{v(2r-1)},\ldots,E_{v\left((m-1)r-1\right)}$
form an edge path $\rho'$ that contains all vertices of $\Gamma$. We can
lift $\rho'$ to a path $\rho$ in $T$ that is contained in the axis of $t$.
We have found a lift of the edges $E_0,\ldots,E_{m-2}$, i.e., we have found
a lift of a spanning tree of $\Gamma$.

\begin{figure}
\renewcommand{\epsfsize}[2]{0.9\textwidth}
\psfrag{e0}{$E_0$}
\psfrag{e2}{$E_2$}
\psfrag{e4}{$E_4$}
\psfrag{e5}{$E_5$}
\psfrag{e7}{$E_7$}
\psfrag{e8}{$E_8$}
\psfrag{dotdot}{$\cdots$}
\psfrag{v8}{$V_8$}
\psfrag{v0}{$V_0$}
\psfrag{v2}{$V_2$}
\psfrag{v3}{$V_3$}
\psfrag{v4}{$V_4$}
\psfrag{v5}{$V_5$}
\psfrag{v7}{$V_7$}
\psfrag{t}{$t$}
\psfrag{a}{$a$}
\psfrag{va}{$w_1a$}
\center{\epsfbox{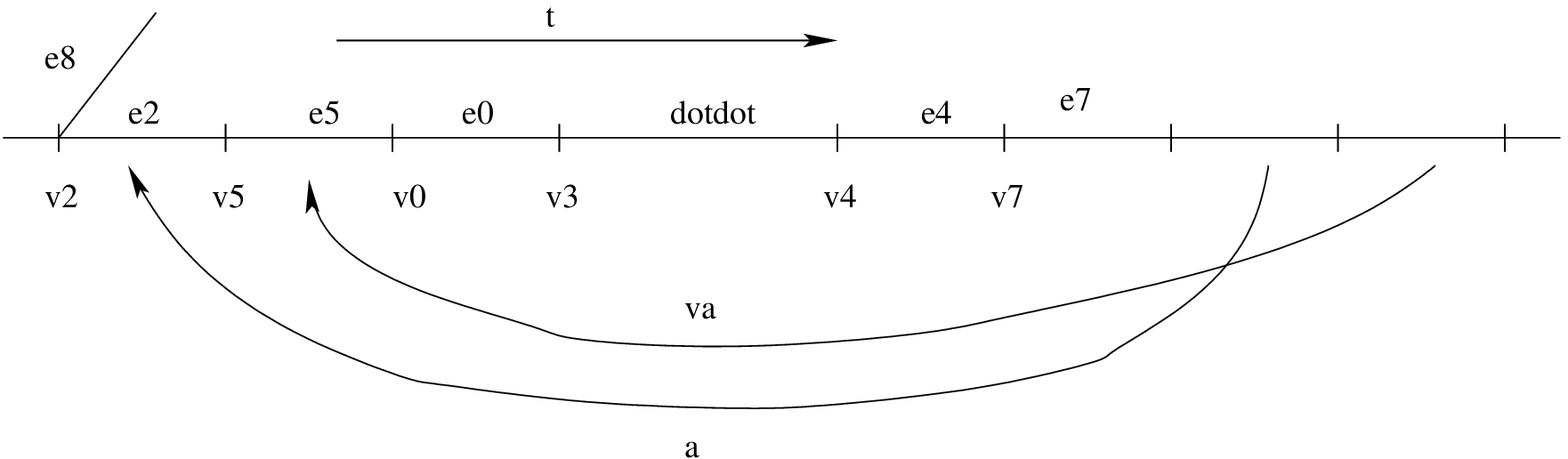}}
\caption{A lift of $\Gamma$ for $m=8, r=s=3$.}\label{hnnhyp2pic}
\end{figure}

The subpath of $t\rho$ that does not overlap with $\rho$ consists of the
edges
$tE_{m-2},tE_{r-2},\ldots,tE_{v\left((s-1)r-2\right)}=E_{m-r-1}$.
Let $E_{m-1}=tE_{m-2}$, and choose $a\in F$ such that $atE_{r-2}=E_{r-1}$.
(\figref{hnnhyp2pic} illustrates these choices in the case $m=8, r=s=3$.)
Note that there exists some $w_1\in F_{\tau(E_{r-1})}$ such that
$w_1atE_{v(2r-2)}=E_{v(2r-1)}$ (\figref{hnnhyp2pic}).
Similarly, for $j=1,\ldots,s-2$, we can
find some $w_j\in F_{\tau(E_{v(jr-1)})}$ such that
$w_j \cdots w_1atE_{v\left((j+1)r-2\right)}=E_{v\left((j+1)r-1\right)}$.

Using the words $a,w_1a,\ldots,w_{s-2}\cdots w_1a$, we now recursively
construct a lift of the remaining edges of $\Gamma$. Then we find
generators of $F$ and their images as before. The details are
left to the reader.

\subsection{Amalgamated free product}

We assume that $M_\phi$ splits as an amalgamated free
product over $\mathbb{Z}$, i.e., $M_\phi=G_1\ast_{\mathbb{Z}}G_2$.
In particular, $t$ acts on $\Gamma$ with two orbits of vertices and
one orbit of edges. As before, we distinguish three subcases depending on
the translation length of the action of $t$ on $T$.

\subsubsection{Elliptic $t$-action}

\begin{figure}
\renewcommand{\epsfsize}[2]{0.8\textwidth}
\psfrag{v}{$V$}
\psfrag{w0}{$W$}
\psfrag{w1}{$tW$}
\psfrag{w2}{$t^2W$}
\psfrag{vv}{$v$}
\psfrag{E}{$E$}
\psfrag{tE}{$tE$}
\psfrag{t2E}{$t^2E$}
\psfrag{t3E}{$t^3E$}
\psfrag{t4E}{$t^4E$}
\psfrag{t5E}{$\cdots$}
\psfrag{t6E}{$t^{kn}E$}
\center{\epsfbox{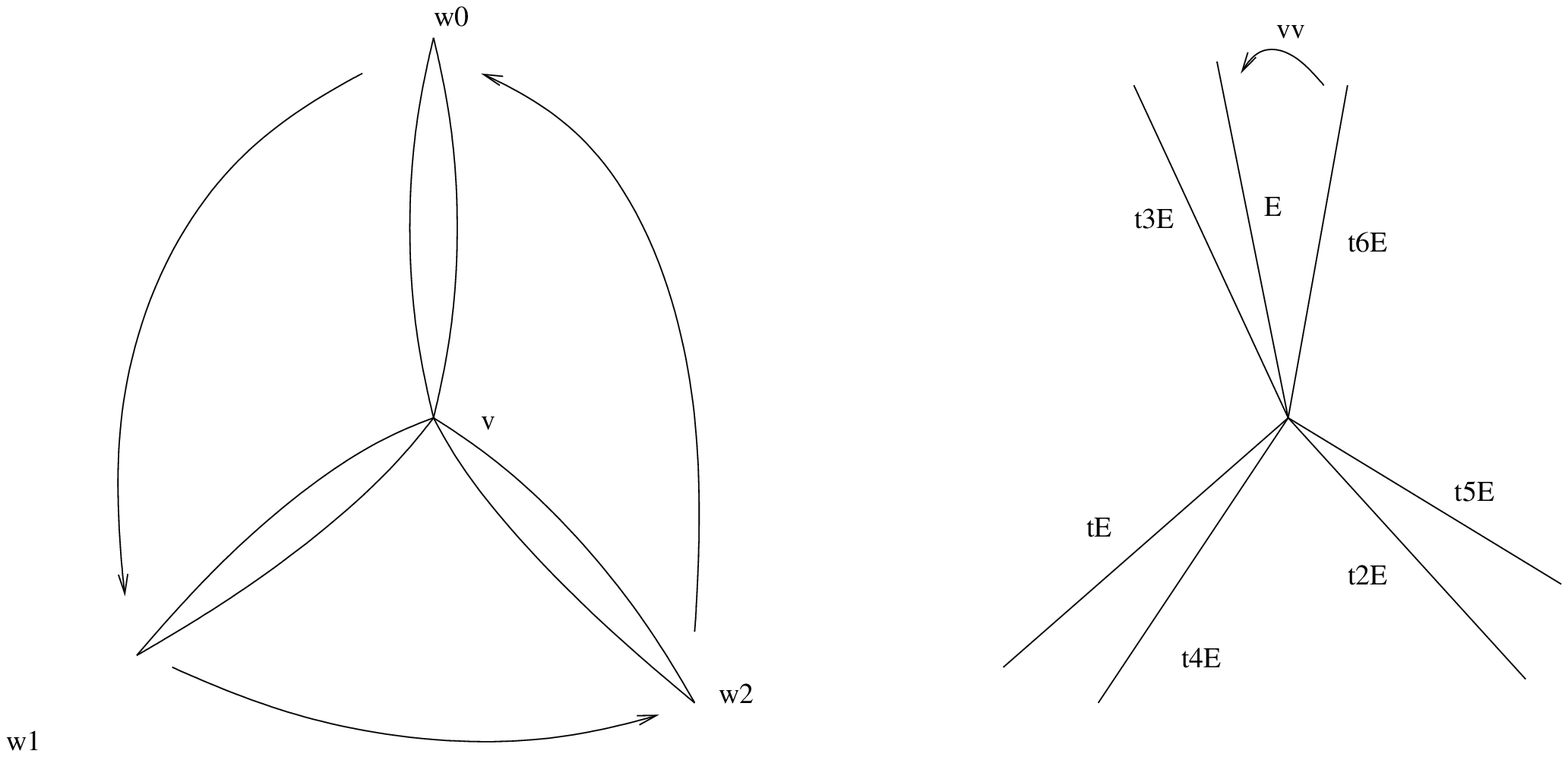}}
\caption{Amalgamated product, elliptic $t$-action.}\label{case3pic}
\end{figure}

We consider graphs $ \Gamma $ with two orbits of $1$, $n$ vertices
respectively (\figref{case3pic}), and one orbit of $kn$ edges.
Then we can express $F$ as
\[
	F=F_V\ast \left(\prod_{i=0}^{n-1} F_{t^iW}\right)\ast
		\langle a_n,\ldots,a_{kn-1} \rangle
\]
such that
\begin{align*}
\phi(F_V)&=F_V\\
\phi(F_{t^iW})&=\begin{Cases}
F_{t^{i+1}W} & if $0\leq i<n-1$;\\
a_n^{-1}F_Wa_n & if $i=n-1$;\\
\end{Cases}\\
\phi(a_i)&=\begin{Cases}
a_{i+1} & if $i<kn-1$;\\
a_{(k-1)n}^{-1}a_{(k-2)n}^{-1}\cdots a_n^{-1}wv & if $i=kn-1$,
where $v\in F_V$ and $w\in F_W$.\\
\end{Cases}
\end{align*}

Note that if $k=1$, then the factor $\langle a_n,\ldots,a_{kn-1} \rangle$
is trivial.

\begin{proof}
As usual, we begin by lifting $\Gamma$ to $T$. Let $V$ be the (unique)
vertex with $tV=V$, and let $E$ be some edge with $\iota(E)=V$. Let
$W=\tau(E)$. Let $E_0=E$ and $E_{i+1}=tE_i$ for $i<kn$. Then
the union $S$ of $E_0,\ldots,E_{kn-1}$ is the desired lift of $\Gamma$.
There exists some $v\in F_V$ such that $vtE_{kn-1}=E_0$.
Choose $a_i\in F$ such that $a_i\tau(E_i)=\tau(E_{i-n})$ for
$n\leq i<kn$.

Using computations analogous to those in previous sections, we can
immediately read off that $\phi(a_i)=a_{i+1}$ if $i<kn$. Moreover,
we have $a_{kn-1}t^{kn-1}W=t^{(k-1)n-1}W$, which implies that
\[
\phi(a_{kn-1})t^{kn}W=t^{(k-1)n}W=a_{(k-1)n}^{-1}t^{(k-2)n}W
=a_{(k-1)n}^{-1}a_{(k-2)n}^{-1}\cdots a_n^{-1}W.
\]
Using the identity $vt^{kn}W=W$, we immediately see that
\[
a_na_{2n}\cdots a_{(k-1)n}\phi(a_{kn-1})v^{-1}=w
\]
for some $w\in F_W$, hence
\[
\phi(a_{kn-1})=a_{(k-1)n}^{-1}a_{(k-2)n}^{-1}\cdots a_n^{-1}wv.
\]
\end{proof}

As before, we find a sequence of Tietze moves that shows that the mapping
torus of an automorphism of this form splits over $\mathbb{Z}$.

\begin{equation*}
\begin{split}
M_\phi &= \langle F_V, F_W, \ldots, F_{t^{n-1}}W, 
a_n, \ldots, a_{kn-1}, t ~\vert~ tF_Vt^{-1}=F_V, \ldots\rangle \\
&\cong \langle F_V, F_W, a_n, \ldots, a_{kn-1}, s, t~\vert~tF_Vt^{-1}=F_V,
s=a_nt^n, sF_Ws^{-1}=F_W,\ldots \rangle. \\
\end{split}
\end{equation*}

The crucial observation at this point is that the relation
$t^{kn}a_nt^{-kn}=\phi^{kn}(a_n)$
can be rewritten as $w^{-1}s^k=vt^{kn}$, using $s=a_nt^n$ and the
structure of $\phi$.  After a few more Tietze moves, we can explicitly
see the splitting.

\begin{equation*}
\begin{split}
M_\phi &\cong \langle F_V, F_W, s, t~\vert~tF_Vt^{-1}=F_V,
sF_Ws^{-1}=F_W, w^{-1}s^k=vt^{kn} \rangle \\
&\cong \langle F_V, t~\vert~tF_Vt^{-1} \rangle
\ast_{\langle vt^{nk}=w^{-1}s^k \rangle}
\langle F_W, s~\vert~sF_Ws^{-1} \rangle. \\
\end{split}
\end{equation*}

\begin{rem}
This case includes the obvious splittings that arise from
automorphisms that respect a free product decomposition of $F$.
\end{rem}

\begin{ex}
This case contains certain {\em inessential} splittings \cite{jsj}, i.e.,
splittings with edge groups that are cyclic, but not maximally cyclic.
Consider, for example, the automorphism $\phi: F(x,y)\rightarrow F(x,y)$,
$x\mapsto y$, $y\mapsto x$. Then $M_\phi$ admits an inessential splitting
over $\mathbb{Z}$:

\begin{equation*}
\begin{split}
M_\phi &= \langle x, y, t~\vert~txt^{-1}=y, tyt^{-1}=x \rangle \\
&\cong \langle x, s, t~\vert~s=t^2, sxs^{-1}=x \rangle \\
&\cong \langle x, s~\vert~sxs^{-1}=x \rangle
\ast_{\langle s\sim t^2\rangle} \langle t \rangle.\\
\end{split}
\end{equation*}

In this example, $F_V$ is trivial, and $\Gamma$ has two edges.
\end{ex}

\subsubsection{Amalgamated free product: Two orbits of edges, special case}

\begin{figure}
\renewcommand{\epsfsize}[2]{0.9\textwidth}
\psfrag{a}{$a$}
\psfrag{V}{$V$}
\psfrag{tV}{$tV$}
\psfrag{tm1V}{$t^{m-1}V$}
\psfrag{tmV}{$t^mV$}
\psfrag{W}{$W$}
\psfrag{tW}{$tW$}
\psfrag{tm1W}{$t^{m-1}W$}
\psfrag{E0}{$E_0$}
\psfrag{En}{$E'$}
\psfrag{E1}{$E_1$}
\psfrag{Em1}{$E_{m-1}$}
\psfrag{Emn1}{$E_{n+m-1}$}
\center{\epsfbox{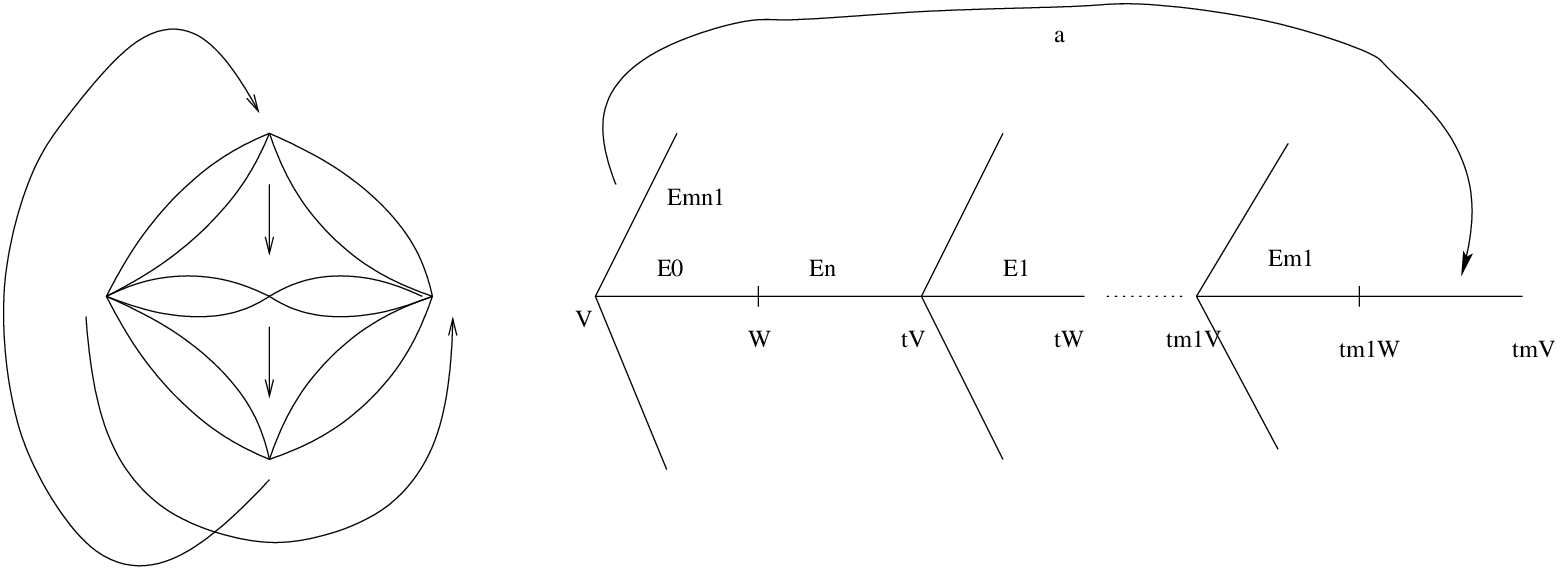}}
\caption{Amalgamated product, hyperbolic $t$-action with translation
distance $2$.}\label{case4pic}
\end{figure}

We consider graphs $\Gamma$ with two orbits of $m$, $n$ vertices respectively,
and $kmn$ edges (\figref{case4pic}). Note that $m$ and $n$ are necessarily
coprime. We may assume $m<n$. As above, we define $v(i)$ to be the remainder
of $i$ under division by $m$, and we let $w(i)$ equal the remainder of $i$
under division by $n$.

We first assume that $n\equiv 1\mod m$. In this case, we have
\[
	F=\left(\prod_{i=0}^{m-1} F_{t^iV}\right)\ast
		\left(\prod_{i=0}^{n-1} F_{t^iW}\right)\ast
		\langle a_{n+m-1},\ldots,a_{kmn-1} \rangle
\]
where
\begin{align*}
\phi(F_{t^iV})&=\begin{Cases}
F_{t^{i+1}V}  & if $i<m-1$;\\
a_{n+m-1}^{-1}F_Va_{n+m-1} & if $i=m-1$;\\
\end{Cases}\\
\phi(F_{t^iW})&=\begin{Cases}
F_{t^{i+1}W} & if $i<n-1$;\\
xF_Wx^{-1} & if $i=n-1$;
\end{Cases}
\end{align*}
for
\[
x=\phi(a^{-1}_{n+m-1})\phi^{m+1}(a^{-1}_{n+m-1})\phi^{2m+1}(a^{-1}_{n+m-1})
\cdots\phi^{n-m}(a^{-1}_{n+m-1})
\]
Moreover, we have
\[
\phi(a_i)=\begin{Cases}
a_{i+1} & if $v(i)\neq m-1$, $v(i-n)\neq m-1$;\\
a_{n+m-1}a_{i+1} & if $v(i-n)=m-1$;\\
a_{i+1}a_{n+m-1}^{-1} & if $v(i)=m-1$, $i\neq knm-1$; \\
a_{(km-1)n}^{-1}a_{(km-2)n}^{-1}\cdots & \\
a_{3n}^{-1}a_{2n}^{-1}wva_{n+m-1}^{-1} & if $i=kmn-1$,\\
& for some $v\in F_V$ and $w\in F_W$.
\end{Cases}
\]

\begin{proof}
As in the previous cases, we first construct a lift of the graph
$\Gamma$. Pick some vertex $V\in T$ and some edge $E$ such that
$\iota(E)=V$. Let $W=\tau(E)$. Since $n\equiv 1\mod m$, there
exists some edge $E'$
such that $\iota(E')=W$ and $E'$ is in the $F$-orbit of
$t^nE$. Let $E_0=E$. There is no loss in assuming
that $tV=\iota(E')$ (\figref{case4pic}).

There exists some $a\in F$ such that $t^mV=aV$. For $0\leq i<knm-1$,
we now recursively define $E_{i+1}$ by letting $E_{i+1}=a^{-1}tE_i$
if $v(i)=m-1$, $E_{i+1}=tE_i$ otherwise. There exists some $v\in F_V$
such that $v(a^{-1}t^m)^{nk}E_0=E_0$. Moreover, there is no loss in
assuming that $t(a^{-1}t^m)^{[\frac{n}{m}]}W=W$, i.e., $E_n=E'$.
(here $[\frac{n}{m}]$ denotes the integral part of $\frac{n}{m}$)

We now choose $a_i$ for $n+m-1\leq i \leq knm-1$ such that
$a_i\tau(E_i)=\tau(E_{i-n})$. In particular, we can let $a_{n+m-1}=a$.
The elements $a_{n+m-1},\ldots,a_{knm-1}$ and the vertex groups
generate $F$. With these definitions, a careful analysis of the action
of the $a_i$s on $T$ yields the above description of $\phi$. The
computations are essentially the same as before.
\end{proof}

In order to see that an automorphism of this form gives rise to a
mapping torus that splits over $\mathbb{Z}$, we first rewrite the
presentation of $M_\phi$ as
\begin{align*}
M_\phi\cong \langle F_V, F_W, a, t~|~& (a^{-1}t^m)F_V(a^{-1}t^m)^{-1}=F_V, \\
& (t(a^{-1}t^m)^{[\frac{n}{m}]})F_W(t(a^{-1}t^m)^{[\frac{n}{m}]})^{-1}=F_W, \\
& t^{knm}at^{-knm}=\phi^{knm}(a) \rangle \\
\end{align*}
Now we introduce the letters $r=a^{-1}t^m$ and
$s=t(a^{-1}t^m)^{[\frac{n}{m}]}$, and we observe that the relation
$t^{knm}at^{-knm}=\phi^{knm}(a)$ can be rewritten as
$vr^{kn}=w^{-1}s^{km}$. Moreover, since $t=sr^{-[\frac{n}{m}]}$
and $a=t^mr^{-1}$, we can eliminate $a$ and $t$ and see the splitting
of $M_\phi$:
\[
M_\phi \cong
\langle F_V, r~|~rF_Vr^{-1}=F_V \rangle
\ast_{\langle vr^{kn}\sim w^{-1}s^{km} \rangle}
\langle F_W, s~|~sF_Ws^{-1}=F_W \rangle.
\]

\subsubsection{Amalgamated free product: Two orbits of edges, general case}

Finally, we need to deal with the case where $n\not\equiv 1\mod m$.
To this end, choose $s$ such that $sn\equiv 1\mod m$ and
$1<s<m$. We construct a lift of $\Gamma$ to $T$ as above,
except we start with a pair of edges $E$ and $E'$ in the
$F$-orbit of $E, t^{sn}E$ respectively. A spanning tree in this case
consists of the edges $E_0,\ldots,E_{n-1}$ and $E_{sn},\ldots,E_{sn+m-2}$.
Hence, we obtain generators $a_n,\ldots,a_{sn-1}$ and $a_{sn+m-1},
\ldots,a_{knm-1}$. Our usual computation of images now gives us the
structure of the automorphism. The details are left to the reader.

\subsection{Summary}

The above discussion covers all possible cases. We have obtained the
main result of this paper.

\begin{thmnp}\label{splitthm}
$M_\phi $ splits over $ \mathbb{Z} $
if and only if $ \phi $ fits into one of the cases listed above.
\end{thmnp}

We conclude this section with a corollary that shows that there
is no shortage of automorphisms whose mapping tori do {\em not}
split over $ \mathbb{Z} $. This corollary was also proved in
\cite{kk}, without a general characterization of splittings
over $\mathbb{Z}$. Recall that an automorphism $\phi: F\rightarrow F$
is called {\em irreducible with irreducible powers} if no positive power
of $\phi$ preserves the conjugacy class of a proper free factor of $F$.

\begin{cor}
If $\phi$ is irreducible with irreducible powers, then $M_\phi$ does
not split over $\mathbb{Z}$.
\end{cor}

%% file: examples.tex
\section{Examples}\label{exsec}

We list some examples of hyperbolic mapping tori that split over
$\mathbb{Z}$. Note, in particular, that
the following proposition implies that there exist reducible
hyperbolic automorphisms. Moreover, there exist hyperbolic
automorphisms that have train track representatives with
polynomially growing strata.

\begin{prop}\label{splitex}
Let $\phi_i$, $i=1,2$, be a hyperbolic automorphism of a free group $F_i$.
\begin{enumerate}
\item The automorphism $\phi=\phi_1\ast\phi_2$ is hyperbolic.
\item Let $w$ be an element of $F_1$ such that
$$w\phi_1(w)\phi_1^2(w)\cdots\phi_1^{k-1}(w)\neq v\phi_1^k(v^{-1})$$
for all $k\geq 1$ and $v\in F_1$. Then the automorphism $\psi$ of
$F=F_1\ast <a>$ defined by $\psi|_{F_1}=\phi_1$, $\psi(a)=aw$, is
hyperbolic. \label{polystrat}
\end{enumerate}
\end{prop}

\begin{proof}
\begin{enumerate}
\item Suppose $\phi=\phi_1\ast\phi_2$ is not hyperbolic.
Since an automorphisms is hyperbolic if and only if it
is atoroidal by \thmref{thesismain}, there exists some
word $1\neq w\in F_1\ast F_2$ such
that $\psi^M(w)$ is conjugate to $w$ for some $M\geq 1$.
There is no loss in assuming that
$w=w_0w_1\cdots w_{k-1}$, where $k$ is even, $1\neq w_{2i}\in F_1$ and
$1\neq w_{2i+1}\in F_2$.

Clearly, both $w$ and $\psi^M(w)$ are cyclically
reduced, and there exists some even number $m$ such that
$\phi^M(w_i)=w_{i+m}$ (indices modulo $k$). But this implies that
        $\psi^{kM}(w_i)=w_i$,
hence neither $\phi_1$ nor $\phi_2$ are atoroidal, so they are not
hyperbolic, which contradicts our hypothesis.

\item Suppose that $\phi$ is not hyperbolic. Then, by \thmref{thesismain},
there exists some word $1\neq w\in F$ such that $\phi^M(w)$ is conjugate
to $w$ for some $M>1$. We can write $w=w_0w_1\cdots w_{k-1}$ where either
$w_i\in F_1$, or $w_i=ax$ for some $x\in F_1$, or $w_i=xa^{-1}$, or
$w=axa^{-1}$. Moreover, we may assume that $w$ is cyclically reduced and
that $k$ is minimal. In particular, there is no cancellation between
successive subwords $w_i$, $w_{i+1}$ (indices modulo $k$).

Then for any $m>1$, there is no cancellation between $\phi^m(w_i)$ and
$\phi^m(w_{i+1})$, and $\phi^m(w)$ is cyclically reduced. Moreover, the
image of $w_i$ is of the same form as $w_i$, e.g.,
a subword of the form $ax$ is mapped to a subword of the form $ax'$.
This implies that $\phi^{kM}$ maps each subword $w_i$ to a conjugate
of itself. Since $\phi$ is hyperbolic, this rules out subwords
$1\neq w_i\in F_1$ as well as subwords
of the form $w_i=axa^{-1}$.

The choice of $w$ also rules out subwords of the form $ax$ and $xa^{-1}$,
which implies $w=1$, a contradiction.

\end{enumerate}
\end{proof}

\begin{ex}
\begin{enumerate}
\item The class of PV-automorphisms introduced in \cite{PV} provides a
rich source of atoroidal automorphisms. For example, for
$F_1=\langle~x,y,z~\rangle$, the automorphism
$\alpha:F_1\rightarrow F_1$ defined by $\alpha(x)=y$, $\alpha(y)=z$, and
$\alpha(z)=xy$ is irreducible and atoroidal.
\item Let $\phi=\alpha^3$. By abelianizing, we can easily see that
$w=x$ satisfies the hypothesis of \propref{splitex}, Part~\ref{polystrat},
so we obtain an explicit example of a hyperbolic automorphism with
a train track representative with a stratum of polynomial growth \cite{hb1}.
\end{enumerate}
\end{ex}